\documentclass[11pt]{article}
\usepackage{mathrsfs}
\usepackage{amssymb,amsfonts}
 \usepackage[bookmarksnumbered, colorlinks, plainpages]{hyperref}
\usepackage{cite}
\def\a{\alpha}
\def\b{\beta}

\def\ci{\circ}

\def\D{\Delta}

\def\g{\gamma}

\def\lr{\longrightarrow}
\def\o{\otimes}

\def\s{\sigma}

\def\v{\varepsilon}

\def\1{^{-1}}
\def\2{^{-2}}
\def\3{^{-3}}
\date{}
\newcounter{zlist}
  
\begin{document}
 \renewcommand{\baselinestretch}{1.2}
 \renewcommand{\arraystretch}{1.0}

 \title{\bf A construction of Hom-Yetter-Drinfeld category}
 \author{{\bf Haiying Li, Tianshui Ma\footnote{Corresponding author}}\\
 {\small  School of Mathematics and Information Science, Henan Normal University,}\\
 {\small Xinxiang 453007, China}\\
 {\small E-mail: haiyingli2012@yahoo.com (H. Li);~~matianshui@yahoo.com (T. Ma)}}
 \maketitle

 \begin{center}
 \begin{minipage}{12.cm}

 {\bf \begin{center}ABSTRACT \end{center}}

 In continuation of our recent work about smash product Hom-Hopf algebras in \cite{MLY}, we introduce Hom-Yetter-Drinfeld category $_H^H{\mathbb{YD}}$ via Radford biproduct Hom-Hopf algebra, and prove that the Hom-Yetter-Drinfeld modules can provide solutions of the Hom-Yang-Baxter equation and $_H^H{\mathbb{YD}}$ is a pre-braided tensor category, where $(H, \b, S)$ is a Hom-Hopf algebra. Furthermore, we obtain that $(A^{\natural}_{\diamond} H,\a\o \b)$ is a Radford biproduct Hom-Hopf algebra if and only if  $(A,\a)$ is a Hopf algebra in the category $_H^H{\mathbb{YD}}$. At last, some examples and applications are given.
 \vskip 0.1cm



 \end{minipage}
 \end{center}
 \normalsize

\section{Introduction }
\def\theequation{1. \arabic{equation}}
\setcounter{equation} {0} \hskip\parindent
 The motivation to introduce Hom-type algebras comes for examples related to q-deformations of Witt and Virasoro algebras, which play an important role in physics, mainly in conformal field theory. Hom-structures (Lie algebras, algebras, coalgebras, Hopf algebras) have been intensively investigated in the literature recently, see \cite{CG,CWZ1,CZ,LM,LS,MLL,MLY,MP1,MP2,MS1,MS2,Yau1,Yau2,Yau3,Yau4,Yau5,Yau6,Yau7,Zhang,Zheng}. Hom-algebras are generalizations of algebras obtained by a twisting map, which have been introduced for the first time in \cite{MS1} by Makhlouf and Silvestrov. The associativity is replaced by Hom-associativity, Hom-coassociativity for a Hom-coalgebra can be considered in a similar way.

 In \cite{Yau1,Yau5}, Yau introduced and characterized the concept of module Hom-algebras as a twisted version of usual module algebras and the dual version (i.e. comodule Hom-coalgebras) was studied by Zhang in \cite{Zhang}. Based on Yau's definition of module Hom-algebras, Ma-Li-Yang in \cite{MLY} constructed smash product Hom-Hopf algebra $(A\natural H,\a\o \b)$ generalizing the Molnar's smash product (see \cite{MW1}), and gave the cobraided structure (in the sense of Yau's definition in \cite{Yau4}) on $(A\natural H,\a\o \b)$, and also considered the case of twist tensor product  Hom-Hopf algebra. Makhlouf and Panaite defined and studied a class of Yetter-Drinfeld modules over Hom-bialgebras in \cite{MP1} and derived the constructions of twistors, pseudotwistors, twisted tensor product and smash product in the setting of Hom-case in \cite{MP2}.

 Yetter-Drinfeld modules are known to be at the origin of a very vast family of solutions to the Yang-Baxter equation.  Let $H$ be a bialgebra, $A$ a left $H$-module algebra and a left $H$-comodule coalgebra. In \cite{Ra1}, Radford gave a construction of bialgebra (called Radford biproduct bialgebra) by combining the smash product algebra $A\# H$ with the smash coproduct coalgebra $A\times H$. Majid (see \cite{Maj1,Maj2}) made the following conclusion: $A$ is a bialgebra in Yetter-Drinfeld category $_H^H{\cal YD}$ if and only if $A\star H$ is a Radford biproduct. The Radford biproduct plays an important role in the lifting method for the classification of finite dimensional pointed Hopf algebras (see \cite{AS}).

 In this paper, we introduce Hom-Yetter-Drinfeld category $_H^H{\mathbb{YD}}$ via Radford biproduct Hom-Hopf algebra, and prove that the Hom-Yetter-Drinfeld modules can provide solutions of the Hom-Yang-Baxter equation. Furthermore, we obtain that $(A^{\natural}_{\diamond} H,\a\o \b)$ is a Radford biproduct Hom-Hopf algebra if and only if  $(A,\a)$ is a Hom-Hopf algebra in the category $_H^H{\mathbb{YD}}$.

 This article is organized as follows. In Section 2, we recall some definitions and results which will be used later. Let $(H, \b)$ be a Hom-bialgebra, $(A, \a)$ a left $(H, \b)$-module Hom-algebra and a left $(H, \b)$-comodule Hom-coalgebra. In \cite{MLY}, the smash product Hom-algebra $(A\natural H,\a\o \b)$ was constructed.  In Section 3, we first define smash coproduct Hom-coalgebra $(A\diamond H,\a\o \b)$ (see Proposition 3.1), then derive necessary and sufficient conditions for $(A\natural H,\a\o \b)$ and $(A\diamond H,\a\o \b)$ to be a Hom-bialgebra, which is called Radford biproduct Hom-bialgebra and denoted by $(A^{\natural}_{\diamond} H,\a\o \b)$ (see Theorem 3.3, 3.6). In Section 4, we introduce the concept of Hom-Yetter-Drinfeld category $_H^H{\mathbb{YD}}$ (see Definition 4.1,4.2), which is different from the one defined by Makhlouf and Panaite in \cite{MP1},  the one defined by Chen and Zhang in \cite{CZ} and the one defined by Liu and Shen in \cite{LS}.  We also prove that the Hom-Yetter-Drinfeld modules can provide solutions of the Hom-Yang-Baxter equation in the sense of Yau's definition in \cite{Yau3,Yau7,Yau6} (see Proposition 4.3) and that $_H^H{\mathbb{YD}}$ is a pre-braided tensor category (see Theorem 4.7). Furthermore, we obtain that $(A^{\natural}_{\diamond} H,\a\o \b)$ is a Radford biproduct Hom-Hopf algebra if and only if  $(A,\a)$ is a Hom-Hopf algebra in the category $_H^H{\mathbb{YD}}$ (see Theorem 4.8), which generalizes the Majid's result in \cite{Maj1,Maj2}.  In last section, some examples and applications are given.

 Throughout this paper we freely use the Hopf algebras and coalgebras terminology introduced in \cite{DNR,Ra2,S,Sw}.

 The authors were informed by the Editor that the following paper \cite{CWZ2} related with  the subject of our paper is accepted for publication.

\section{Preliminaries}
\def\theequation{2. \arabic{equation}}
\setcounter{equation} {0} \hskip\parindent
 Throughout this paper, we follow the definitions and terminologies in \cite{Ka,MLY,Yau1,Yau3,Zhang}, with all algebraic systems supposed to be over the field $K$. Given a $K$-space $M$, we write $id_M$ for the identity map on $M$.
\smallskip

 We now recall some useful definitions.

\smallskip

 {\bf Definition 2.1} A Hom-algebra is a quadruple $(A,\mu,1_A,\a)$ (abbr. $(A,\a)$), where $A$ is a $K$-linear space, $\mu: A\o A \lr A$ is a $K$-linear map, $1_A \in A$ and $\a$ is an automorphism of $A$, such that
 \begin{eqnarray*}
 &(HA1)& \a(aa')=\a(a)\a(a');~~\a(1_A)=1_A,\\
 &(HA2)& \a(a)(a'a'')=(aa')\a(a'');~~a1_A=1_Aa=\a(a)
 \end{eqnarray*}
 are satisfied for $a, a', a''\in A$. Here we use the notation $\mu(a\o a')=aa'$.

 Let $(A, \a)$ and $(B, \b)$ be two Hom-algebras. Then $(A\o B, \a\o \b)$ is a Hom-algebra (called tensor product Hom-algebra) with the multiplication $(a\o b)(a'\o b')=aa'\o bb'$ and unit $1_A\o 1_B$.

\smallskip

 {\bf Definition 2.2} A Hom-coalgebra is a quadruple $(C,\D,\v_C,\b)$ (abbr. $(C,\b)$), where $C$ is a $K$-linear space, $\D: C \lr C\o C$, $\v_C: C\lr K$ are $K$-linear maps, and $\b$ is an automorphism of $C$, such that
 \begin{eqnarray*}
 &(HC1)& \b(c)_1\o \b(c)_2=\b(c_1)\o \b(c_2);~~\v_C\ci \b=\v_C\\
 &(HC2)& \b(c_{1})\o c_{21}\o c_{22}=c_{11}\o c_{12}\o \b(c_{2});~~\v_C(c_1)c_2=c_1\v_C(c_2)=\b(c)
 \end{eqnarray*}
 are satisfied for $c \in A$. Here we use the notation $\D(c)=c_1\o c_2$ (summation implicitly understood).

 Let $(C, \a)$ and $(D, \b)$ be two Hom-coalgebras. Then $(C\o D, \a\o \b)$ is a Hom-coalgebra (called tensor product Hom-coalgebra) with the comultiplication $\D(c\o d)=c_1\o d_1\o c_2\o d_2$ and counit $\v_C\o \v_D$.

\smallskip

 {\bf Definition 2.3} A Hom-bialgebra is a sextuple $(H,\mu,1_H,\D,\v,\g)$ (abbr. $(H,\g)$), where
 $(H,\mu,1_H,\g)$ is a Hom-algebra and $(H,\D,\v,\g)$ is a Hom-coalgebra, such that
 $\D$ and $\v$ are morphisms of Hom-algebras, i.e.
 $$
 \D(hh')=\D(h)\D(h');~~\D(1_H)=1_H\o 1_H,
 $$
 $$
 \v(hh')=\v(h)\v(h');~~\v(1_H)=1.
 $$
 Furthermore, if there exists a linear map $S: H\lr H$ such that
 $$
 S(h_1)h_2=h_1S(h_2)=\v(h)1_H~\hbox{and}~S(\g(h))=\g(S(h)),
 $$
 then we call $(H,\mu,1_H,\D,\v,\g,S)$(abbr. $(H,\g,S)$) a Hom-Hopf algebra.

 Let $(H,\g)$ and $(H',\g')$ be two Hom-bialgebras. The linear map $f: H\lr H'$ is called a Hom-bialgebra map if $f\ci \g=\g'\ci f$
 and at the same time $f$ is a bialgebra map in the usual sense.

 \smallskip

 {\bf Definition 2.4} (see \cite{Yau1,Yau5}) Let $(A,\b)$ be a Hom-algebra. A left $(A,\b)$-Hom-module is a triple $(M,\rhd,\a)$, where $M$ is a linear space, $\rhd: A\o M \lr M$ is a linear map, and $\a$ is an automorphism of $M$, such that
  \begin{eqnarray*}
 &(HM1)& \a(a\rhd m)=\b(a)\rhd \a(m),\\
 &(HM2)& \b(a)\rhd (a'\rhd m)=(aa')\rhd \a(m);~~ 1_A\rhd m=\a(m)
 \end{eqnarray*}
 are satisfied for $a, a' \in A$ and $m\in M$.

 Let $(M,\rhd_M,\a_M)$ and $(N,\rhd_N,\a_N)$ be two left $(A,\b)$-Hom-modules. Then a linear morphism $f: M\lr N$ is called a morphism of left $(A,\b)$-Hom-modules if $f(h\rhd_M m)=h\rhd_N f(m)$ and $\a_M\ci f=f\ci \a_N$.
 \smallskip

 {\bf Remarks} (1) It is obvious that $(A,\mu,\b)$ is a left
 $(A,\b)$-Hom-module.\newline
 \indent{\phantom{\bf Remarks}} (2) When $\b=id_A$ and $\a=id_M$, a left $(A,\b)$-Hom-module
 is the usual left $A$-module.

 \smallskip

 {\bf Definition 2.5} (see \cite{Yau1,Yau5}) Let $(H,\b)$ be a Hom-bialgebra and $(A,\a)$ a Hom-algebra. If $(A,\rhd,\a)$ is a left $(H,\b)$-Hom-module and for all $h\in H$ and $a, a'\in A$,
 \begin{eqnarray*}
 &(HMA1)& \b^{2}(h)\rhd (aa')=(h_1\rhd a)(h_2\rhd a'),\\
 &(HMA2)& h\rhd 1_A=\v_H(h)1_A,
 \end{eqnarray*}
 then $(A,\rhd,\a)$ is called an $(H,\b)$-module Hom-algebra.

 \smallskip

 {\bf Remarks} (1) When $\a=id_A$ and $\b=id_H$, an $(H,\b)$-module Hom-algebra
 is the usual $H$-module algebra.\newline
 \indent{\phantom{\bf Remarks}} (2) Similar to the case of Hopf
 algebras, in \cite{Yau1,Yau5}, Yau concluded that the Eq.$(HMA1)$ is satisfied if and only if $\mu_A$ is a
 morphism of $H$-modules for suitable $H$-module structures on $A\o A$ and $A$, respectively.\newline
 \indent{\phantom{\bf Remarks}} (3)  The smash product Hom-Hopf algebra $(A\natural H,\a\o \b)$ is different from the one defined by Chen, Wang and Zhang in \cite{CWZ1}, since here the construction of $(A\natural B,\a\o \b)$ is based on the concept of the module Hom-algebra introduced by Yau in \cite{Yau1,Yau5}, while two of conditions $(6.1),(6.2)$ in the module Hom-algebra in \cite{CWZ1} are same to the case of Hopf algebra.

 \smallskip

 {\bf Definition 2.6} (see \cite{Zhang}) Let $(C,\b)$ be a Hom-coalgebra. A left $(C,\b)$-Hom-comodule is a triple $(M,\rho,\a)$, where $M$ is a linear space, $\rho:  M\lr C\o M$ (write $\rho(m)=m_{-1}\o m_0,~\forall m\in M$) is a linear map,
 and $\a$ is an automorphism of $M$, such that
 \begin{eqnarray*}
 &(HCM1)&\a(m)_{-1}\o \a(m)_{0}=\b(m_{-1})\o \a(m_{0}),\\
 &(HCM2)&\b(m_{-1})\o m_{0-1}\o m_{00}= m_{-11}\o m_{-12}\o \a(m_{0});~~ \v_C(m_{-1})m_{0}=\a(m)
 \end{eqnarray*}
 are satisfied for all $m\in M$.

 Let $(M,\rho^M,\a_M)$ and $(N,\rho^N,\a_N)$ be two left $(C,\b)$-Hom-comodules. Then a linear map $f: M\lr N$ is called a map of left $(C,\b)$-Hom-comodules if $f(m)_{-1}\o f(m)_{0}=m_{-1}\o f(m_{0})$ and $\a_M\ci f=f\ci \a_N$.

 \smallskip

 {\bf Remarks} (1) It is obvious that $(C,\D_C,\b)$ is a left $(C,\b)$-Hom-comodule.\newline
 \indent{\phantom{\bf Remarks}} (2) When $\b=id_A$ and $\a=id_M$, a left $(C,\b)$-Hom-comodule
 is the usual left $C$-comodule.

 \smallskip

 {\bf Definition 2.7} (see \cite{Zhang}) Let $(H,\b)$ be a Hom-bialgebra and $(C,\a)$ a Hom-coalgebra. If $(C,\rho,\a)$ is a left $(H,\b)$-Hom-comodule and for all $c\in C$,
 \begin{eqnarray*}
 &(HCMC1)&\b^{2}(c_{-1})\o c_{01}\o c_{02}=c_{1-1}c_{2-1}\o c_{10}\o c_{20},\\
 &(HCMC2)&c_{-1}\v_C(c_0)=1_H\v_C(c),
 \end{eqnarray*}
 then $(C,\rho,\a)$ is called an $(H,\b)$-comodule Hom-coalgebra.

 \smallskip

 {\bf Remarks} (1) When $\a=id_A$ and $\b=id_H$, an $(H,\b)$-comodule Hom-coalgebra
 is the usual $H$-comodule coalgebra.\newline
 \indent{\phantom{\bf Remarks}} (2) Similar to the case of Hopf
 algebras, in \cite{Zhang}, Zhang concluded that the Eq.$(HCMC1)$ is satisfied if and only if $\D_C$ is a
 morphism of $H$-comodules for suitable $H$-comodule structures on $C\o C$ and $C$, respectively.

 \smallskip

 {\bf Definition 2.8} (see \cite{MLY}) Let $(H,\b)$ be a Hom-bialgebra and $(C,\a)$ a Hom-coalgebra. If $(C,\rhd,\a)$ is a left $(H,\b)$-Hom-module and for all $h\in H$ and $c\in A$,
 \begin{eqnarray*}
 &(HMC1)&(h\rhd c)_1\o (h\rhd c)_2=(h_1\rhd c_1)\o (h_2\rhd c_2),\\
 &(HMC2)&\v_C(h\rhd c)=\v_H(h)\v_C(c),
 \end{eqnarray*}
 then $(C,\rhd,\a)$ is called an $(H,\b)$-module Hom-coalgebra.

 \smallskip

 {\bf Remark} When $\a=id_C$ and $\b=id_H$, an $(H,\b)$-module Hom-coalgebra is the usual $H$-module coalgebra.

 \smallskip

 {\bf Definition 2.9} (see \cite{Yau2}) Let $(H,\b)$ be a Hom-bialgebra and $(A,\a)$ a Hom-algebra. If $(A,\rho,\a)$ is a left $(H,\b)$-Hom-comodule and for all $a, a'\in A$,
 \begin{eqnarray*}
 &(HCMA1)&\rho(aa')=a_{-1}a'_{-1}\o a_{0}a'_{0},\\
 &(HCMA2)&\rho(1_A)=1_H\o 1_A,
 \end{eqnarray*}
 then $(A,\rho,\a)$ is called an $(H,\b)$-comodule Hom-algebra.

 \smallskip

 {\bf Remark} When $\a=id_A$ and $\b=id_H$, an $(H,\b)$-comodule Hom-algebra is the usual $H$-comodule algebra.

 \smallskip

 {\bf Definition 2.10} (see \cite{MLY}) Let $(H,\b)$ be a Hom-bialgebra and $(A,\rhd,\a)$ an $(H,\b)$-module Hom-algebra. Then $(A\natural H, \a\o \b)$ ($A\natural H=A\o H$ as a linear space) with the multiplication
 $$
 (a\o h)(a'\o h')=a(h_1\rhd \a^{-1}(a'))\o \b^{-1}(h_{2})h',
 $$
 where $a, a'\in A, h, h'\in H$, and unit $1_A\o 1_H$ is a Hom-algebra, we call it smash product Hom-algebra denoted by $(A\natural H,\a\o \b)$.

 \smallskip

 {\bf Remarks} Here the multiplication of smash product Hom-algebra is different from the one defined by Makhlouf and Panaite in Theorem3.1 in \cite{MP2}.

 \smallskip

 {\bf Definition 2.11} (see \cite{AS,Maj2,MP1}) Let $H$ be a bialgebra and $M$ a linear space which is a left $H$-module with action
 $\rhd: H\o M\lr M, h\o m \mapsto h\rhd m$ and a left $H$-comodule with coaction $\rho: M\lr H\o M, \rho(m)=m_{-1}\o m_0$.  Then $M$ is called a (left-left) Yetter-Drinfeld module over $H$ if the following compatibility condition holds, for all $h\in H$ and $m \in M$,
 $$
  (YD) ~~~~~h_1m_{-1}\o (h_2\rhd m_0)=(h_1\rhd m)_{-1}h_2\o (h_1\rhd m)_{0}.
 $$
 When $H$ is a Hopf algebra, then the condition $(YD)$ is equivalent to
 $$
 (YD)' ~~~~~h_1m_{-1}S_H(h_3)\o (h_2\rhd m_0)=(h\rhd m)_{-1} \o (h\rhd m)_{0}.
 $$

 \smallskip

 \section{Radford biproduct Hom-Hopf algebra}
\def\theequation{3. \arabic{equation}}
\setcounter{equation} {0} \hskip\parindent
 In this section, we mainly generalize the Radford biproduct bialgebra in \cite[Theorem1]{Ra1} to the Hom-setting.

 Dual to the Definition 2.10, we have:

 {\bf Proposition 3.1} Let $(H,\b)$ be a Hom-bialgebra and $(C,\rho,\a)$ an $(H,\b)$-comodule Hom-coalgebra. Then $(C\diamond H, \a\o \b)$ ($C\diamond H=C\o H$ as a linear space) with the comultiplication
 $$
 \D_{C\diamond H}(c\o h)=c_1\o c_{2-1}\b^{-1}(h_1)\o \a^{-1}(c_{20})\o h_2,
 $$
 where $c\in C, h\in H$, and counit $\v_C\o \v_H$ is a Hom-coalgebra, we call it smash coproduct Hom-coalgebra denoted by $(C\diamond H,\a\o \b)$.

 In fact, dual to Theorem 3.1 in \cite{MLY}, we have

 {\bf Proposition 3.2} Let $(C,\D_C,\v_C,\a)$ and $(H,\D_H,\v_H,\b)$ be two Hom-coalgebras, $T: C\o H \lr H\o C$ (write $T(c\o h)=h_T\o c_T, \forall c\in C, h\in H$) a linear map such that for all $c\in C, h\in H$,
 $$
 \a(c)_T\o \b(h)_T=\a(c_T)\o \b(h_T).
 $$
 Then $(C\diamond_T H, \a\o \b)$ ($C\diamond_T H=C\o H$ as a linear space) with the comultiplication
 $$
 \D_{C\diamond_T H}(c\o h)=c_1\o \b^{-1}(h_1)_T\o \a^{-1}(c_{2T})\o h_2,
 $$
 and counit $\v_C\o \v_H$ becomes a Hom-coalgebra if and only if the following conditions hold:
 \begin{eqnarray*}
 &(C1)& \v_H(h_T)c_T=\v_H(h)\a(c);~~h_T\v_C(c_T)=\b(h)\v_C(c),\\
 &(C2)& h_{T1}\o h_{T2}\o \a(c_T)=\b(\b^{-1}(h_1)_T)\o h_{2t}\o c_{Tt},\\
 &(C3)& \b(h_T)\o \a(c)_{T1}\o \a(c)_{T2}=h_{Tt}\o \a(c_1)_t\o \a(c_{2T}),
 \end{eqnarray*}
 where $c\in C, h\in H$ and $t$ is a copy of $T$.

 We call this Hom-coalgebra $T$-smash coproduct Hom-coalgebra and denote it by $(C\diamond_T H,\a\o \b)$.

 \smallskip

 {\bf Remarks} (1) Let $T(c\o h)=c_{-1}h\o c_0$ in $C\diamond_T H$, we can get the smash coproduct Hom-coalgebra $C\diamond H$. \newline
 \indent{\phantom{\bf Remarks}}(2) Here the comultiplication of $T$-smash coproduct Hom-coalgebra is slightly different from the one defined by Zheng in \cite{Zheng}. And the conditions $(C1)-(C3)$ are simpler than the ones in \cite{Zheng}.

 \smallskip

 {\bf Theorem 3.3}  Let $(H, \b)$ be a Hom-bialgebra, $(A, \a)$ a left $(H, \b)$-module Hom-algebra with module structure $\rhd: H\o A\lr A$ and a left $(H, \b)$-comodule Hom-coalgebra with comodule structure $\rho: A\lr H\o A$. Then the following are equivalent:
 \begin{itemize}
 \item $(A^{\natural}_{\diamond} H, \mu_{A\natural H}, 1_A\o 1_H, \D_{A\diamond H}, \v_A\o \v_H, \a\o \b)$ is a Hom-bialgebra, where
 $(A\natural H, \a\o \b)$ is a smash product Hom-algebra and $(A\diamond H, \a\o \b)$ is a smash coproduct Hom-coalgebra.

 \item The following conditions hold ($\forall~a,b\in A$ and $h\in H$):

 \quad (R1) $(A,\rho,\a)$ is an $(H,\b)$-comodule Hom-algebra,

 \quad (R2) $(A,\rhd,\a)$ is an $(H,\b)$-module Hom-coalgebra,

 \quad (R3) $\v_A$ is a Hom-algebra map and $\D_A(1_A)=1_A\o 1_A$,

 \quad (R4) $\D_A(ab)=a_1(\b^2(a_{2-1})\rhd \a^{-1}(b_1))\o \a^{-1}(a_{20})b_2$,

 \quad (R5) $h_1\b(a_{-1})\o (\b^3(h_2)\rhd a_0)=(\b^2(h_1)\rhd a)_{-1}h_2\o (\b^2(h_1)\rhd a)_{0}$.

 \end{itemize}

 In this case, we call this Hom-bialgebra  Radford biproduct Hom-bialgebra and denote it by $(A^{\natural}_{\diamond} H,\a\o \b)$.

 {\bf Proof} ($\Longleftarrow$)  It is easy to prove that $\v_{A^{\natural}_{\diamond} H}=\v_A\o \v_H$ is a morphism of Hom-algebras. Next we check $\D_{A^{\natural}_{\diamond} H}=\D_{A\diamond H}$  is a morphism of Hom-algebras as follows. For all $a,b\in A$ and $h,g\in H$, we have
 \begin{eqnarray*}
 &&\D_{A^{\natural}_{\diamond} H}((a\o h)(b\o g))\\
 &&~~~~~\stackrel{}{=}(a(h_1\rhd \a\1(b)))_1\o (a(h_1\rhd \a\1(b)))_{2-1}\b\1((\b\1(h_2)g)_{1})\\
 &&~~~~~~~~~~\o \a\1((a(h_1\rhd \a\1(b)))_{20})\o (\b\1(h_2)g)_{2}\\
 &&\stackrel{(HA1)(HC1)}{=}(a(h_1\rhd \a\1(b)))_1\o (a(h_1\rhd \a\1(b)))_{2-1}(\b^{-2}(h_{21})\b\1(g_{1}))\\\
 &&~~~~~~~~~~\o \a\1((a(h_1\rhd \a\1(b)))_{20})\o \b\1(h_{22})g_{2}\\
 &&~~\stackrel{(R4)}{=}a_1(\b^{2}(a_{2-1})\rhd \a\1((h_1\rhd \a\1(b))_{1}))\o (\a\1(a_{20})(h_1\rhd \a\1(b))_{2})_{-1}\\
 &&~~~~~~~~~~\times (\b^{-2}(h_{21})\b\1(g_{1}))\o \a\1((\a\1(a_{20})(h_1\rhd \a\1(b))_{2})_{0})\o \b\1(h_{22})g_{2}\\
 &&~~\stackrel{(HCA1)}{=}a_1(\b^{2}(a_{2-1})\rhd \a\1((h_1\rhd \a\1(b))_{1}))\o (\a\1(a_{20})_{-1}(h_1\rhd \a\1(b))_{2-1})\\
 &&~~~~~~~~~~\times (\b^{-2}(h_{21})\b\1(g_{1}))\o \a\1(\a\1(a_{20})_{0})\a\1((h_1\rhd \a\1(b))_{20})\o \b\1(h_{22})g_{2}\\
 &&~~\stackrel{(HMC1)}{=}a_1(\b^{2}(a_{2-1})\rhd \a\1(h_{11}\rhd \a\1(b_{1})))\o (\a\1(a_{20})_{-1}(h_{12}\rhd \a\1(b_2))_{-1})\\
 &&~~~~~~~~~~\times (\b^{-2}(h_{21})\b\1(g_{1}))\o \a\1(\a\1(a_{20})_{0})\a\1((h_{12}\rhd \a\1(b_2))_{0})\o \b\1(h_{22})g_{2}\\
 &&~~\stackrel{(HA2)}{=}a_1(\b^{2}(a_{2-1})\rhd \a\1(h_{11}\rhd \a\1(b_{1})))\o (\a\1(a_{20})_{-1}\b\1((h_{12}\rhd \a\1(b_2))_{-1}\\
 &&~~~~~~~~~~\times (\b^{-2}(h_{21})))g_{1}\o \a\1(\a\1(a_{20})_{0})\a\1((h_{12}\rhd \a\1(b_2))_{0})\o \b\1(h_{22})g_{2}\\
 &&~~\stackrel{(HC2)}{=}a_1(\b^{2}(a_{2-1})\rhd \a\1(\b(h_{1})\rhd \a\1(b_{1})))\o (\a\1(a_{20})_{-1}\\
 &&~~~~~~~~~~\times \b\1((\b\1(h_{211})\rhd \a\1(b_2))_{-1}\b^{-3}(h_{212})))g_{1}\o \a\1(\a\1(a_{20})_{0})\\
 &&~~~~~~~~~~\times \a\1((\b\1(h_{211})\rhd \a\1(b_2))_{0})\o \b\1(h_{22})g_{2}\\
 &&~~\stackrel{(HC1)}{=}a_1(\b^{2}(a_{2-1})\rhd \a\1(\b(h_{1})\rhd \a\1(b_{1})))\o (\a\1(a_{20})_{-1}\\
 &&~~~~~~~~~~\times \b\1((\b^{2}(\b^{-3}(h_{21})_{1})\rhd \a\1(b_2))_{-1}\b^{-3}(h_{21})_{2}))g_{1}\o \a\1(\a\1(a_{20})_{0})\\
 &&~~~~~~~~~~\times \a\1((\b^{2}(\b^{-3}(h_{21})_{1})\rhd \a\1(b_2))_{0})\o \b\1(h_{22})g_{2}\\
 &&~~\stackrel{(R5)}{=}a_1(\b^{2}(a_{2-1})\rhd \a\1(\b(h_{1})\rhd \a\1(b_{1})))\o (\a\1(a_{20})_{-1}\b\1(\b^{-3}(h_{21})_{1})\\
 &&~~~~~~~~~~\times \b(\a\1(b_2)_{-1})))g_{1}\o \a\1(\a\1(a_{20})_{0}) \a\1(\b^{3}(\b^{-3}(h_{21})_{2})\rhd \a\1(b_2)_{0})\\
 &&~~\stackrel{(HCM1)(HC1)}{=}a_1(\b^{2}(a_{2-1})\rhd \a\1(\b(h_{1})\rhd \a\1(b_{1})))\o (\b\1(a_{20-1})\b\1(\b^{-3}(h_{211})\\
 &&~~~~~~~~~~\times b_{2-1}))g_{1}\o \a^{-2}(a_{200})\a\1(h_{212}\rhd \a\1(b_{20}))\o \b\1(h_{22})g_{2}\\
 &&~~\stackrel{(HCM2)}{=}a_1(\b(a_{2-11})\rhd \a\1(\b(h_{1})\rhd \a\1(b_{1})))\o (\b\1(a_{2-12})\b\1(\b^{-3}(h_{211})\\
 &&~~~~~~~~~~\times b_{2-1}))g_{1}\o \a^{-1}(a_{20})\a\1(h_{212}\rhd \a\1(b_{20}))\o \b\1(h_{22})g_{2}\\
 &&~~\stackrel{(HA2)}{=}a_1(\b(a_{2-11})\rhd \a\1(\b(h_{1})\rhd \a\1(b_{1})))\o (\b\1(a_{2-12})\b^{-3}(h_{211}))\\
 &&~~~~~~~~~~\times (b_{2-1}\b\1(g_{1}))\o \a^{-1}(a_{20})\a\1(h_{212}\rhd \a\1(b_{20}))\o \b\1(h_{22})g_{2}\\
 &&~~\stackrel{(HC2)}{=}a_1(\b(a_{2-11})\rhd \a\1(h_{11}\rhd \a\1(b_{1})))\o (\b\1(a_{2-12})\b^{-2}(h_{12}))(b_{2-1}\b\1(g_{1}))\\
 &&~~~~~~~~~~\o \a^{-1}(a_{20})\a\1(\b(h_{21})\rhd \a\1(b_{20}))\o \b\1(h_{22})g_{2}\\
 &&~~\stackrel{(HM1)}{=}a_1(\b(a_{2-11})\rhd (\b\1(h_{11})\rhd \a^{-2}(b_{1})))\o (\b\1(a_{2-12})\b^{-2}(h_{12}))(b_{2-1}\b\1(g_{1}))\\
 &&~~~~~~~~~~\o \a^{-1}(a_{20})(h_{21}\rhd \a^{-2}(b_{20}))\o \b\1(h_{22})g_{2}\\
 &&~~\stackrel{(HM2)}{=}a_1((a_{2-11}\b\1(h_{11}))\rhd \a\1(b_{1}))\o (\b\1(a_{2-12})\b^{-2}(h_{12}))(b_{2-1}\b\1(g_{1}))\\
 &&~~~~~~~~~~\o \a^{-1}(a_{20})(h_{21}\rhd \a^{-2}(b_{20}))\o \b\1(h_{22})g_{2}\\
 &&~~\stackrel{(HA1)}{=}a_1((a_{2-1}\b\1(h_{1}))_1\rhd \a\1(b_{1}))\o \b\1((a_{2-1}\b\1(h_{1}))_2) (b_{2-1}\b\1(g_{1}))\\
 &&~~~~~~~~~~\o \a^{-1}(a_{20})(h_{21}\rhd \a^{-2}(b_{20}))\o \b\1(h_{22})g_{2}\\
 &&~~\stackrel{}{=}(a_1\o a_{2-1}\b^{-1}(h_1)\o \a^{-1}(a_{20})\o h_2)(b_1\o b_{2-1}\b^{-1}(h_1)\o \a^{-1}(b_{20})\o h_2)\\
 &&~~\stackrel{}{=}\D_{A^{\natural}_{\diamond} H}(a\o h)\D_{A^{\natural}_{\diamond} H}(b\o g),
 \end{eqnarray*}
 and $\D_{A^{\natural}_{\diamond} H}(1_A\o 1_H)=1_A\o 1_H\o 1_A\o 1_H$ can be proved directly.

 $(\Longrightarrow)$  We only verify that the conditions $(R4)$ and $(R5)$ hold, and others hold similarly.
 Since $\D_{A^{\natural}_{\diamond} H}=\D_{A\diamond H}$  is a morphism of Hom-algebras, for all $a,b\in A$ and $h,g\in H$, we have
 \begin{eqnarray*}
 &&a_1((a_{2-1}\b\1(h_{1}))_1\rhd \a\1(b_{1}))\o \b\1((a_{2-1}\b\1(h_{1}))_2) (b_{2-1}\b\1(g_{1}))\\
 &&~~~\o \a^{-1}(a_{20})(h_{21}\rhd \a^{-2}(b_{20}))\o \b\1(h_{22})g_{2}=(a(h_1\rhd \a\1(b)))_1\o (a(h_1\rhd \a\1(b)))_{2-1}\\
 &&~~~\times \b\1((\b\1(h_2)g)_{1})\o \a\1((a(h_1\rhd \a\1(b)))_{20})\o (\b\1(h_2)g)_{2}
 \end{eqnarray*}
 Then, apply $id_A\o \v_H\o id_A\o \v_H$ to the above equation and set $h=g=1_H$, we get $(HB)$. $(HYD)$ can be obtained by using $\v_A\o id_H\o id_A\o \v_H$ to the above equation and setting $a=1_A, g=1_H$.            \hfill $\square$

 \smallskip

 {\bf Remarks} (1) If $\a=id_A$ and $\b=id_H$, then we can get the well-known Radford biproduct bialgebra in \cite[Theorem 1]{Ra1}.\newline
 \indent{\phantom{\bf Remarks}}(2) Theorem 3.3 is different from the one defined by Liu and Shen in \cite{LS}, because the Hom-smash product there is based on the concept of module Hom-algebra in \cite{CWZ1} and ours is based on the Yau's in \cite{Yau1,Yau5}.

 \smallskip

  {\bf Corollary 3.4} (see \cite{MLY}) Let $(A,\a), (H,\b)$ be two Hom-bialgebras, and $(A,\rhd,\a)$ an $(H,\b)$-module Hom-algebra. Then the smash product Hom-algebra $(A\natural H, \a\o \b)$ endowed with the tensor product Hom-coalgebra structure becomes a Hom-bialgebra if and only if $(A,\rhd,\a)$ is an $(H,\b)$-module Hom-coalgebra  and
 $$
 h_1\o h_2\rhd a=h_2\o h_1\rhd a.
 $$

 {\bf Proof} Let the comodule action $\rho$ be trivial, i.e. $\rho(a)=1_H\o \a(a)$ in Theorem 3.3.         \hfill $\square$

 \smallskip

  {\bf Corollary 3.5} Let $(C,\a), (H,\b)$ be two Hom-bialgebras, and $(C,\rho,\a)$ an $(H,\b)$-comodule Hom-coalgebra. Then the smash coproduct Hom-coalgebra $(C\diamond H, \a\o \b)$ endowed with the tensor product Hom-algebra structure becomes a Hom-bialgebra if and only if $(C,\rho,\a)$ is an $(H,\b)$-comodule Hom-algebra  and
 $$
 hc_{-1}\o c_0=c_{-1}h\o c_0.
 $$

 {\bf Proof} Let the module action $\rhd$ be trivial, i.e. $h\rhd c=\v_H(h)\a(c)$ in Theorem 3.3.         \hfill $\square$

 \smallskip

  {\bf Theorem 3.6}  Let $(H,\b,S_{H})$ be a Hom-Hopf algebra, and $(A,\a)$ be a Hom-algebra and a Hom-coalgebra. Assume that $(A^{\natural}_{\diamond} H, \a\o \b)$ is a Radford biproduct Hom-bialgebra defined as above, and $S_A: A\lr A$ is a linear map such that
  $S_A(a_1)a_2=a_1S_A(a_2)=\v_A(a)1_A$ and $\a\ci S_A=S_A\ci \a$ hold.  Then $(A^{\natural}_{\diamond} H, \a\o \b, S_{A^{\natural}_{\diamond} H})$ is a Hom-Hopf algebra, where
  $$
  S_{A^{\natural}_{\diamond} H}(a\o h)=(S_H(a_{-1}\b\1(h))_{1}\rhd S_A(\a\2(a_{0})))\o \b\1(S_H(a_{-1}\b\1(h))_{2}).
  $$

  {\bf Proof} We can compute that $(A^{\natural}_{\diamond} H, \a\o \b, S_{A^{\natural}_{\diamond} H})$ is a Hom-Hopf algebra as follows. For all $a\in A$ and $h\in H$, we have
 \begin{eqnarray*}
 &&(S_{A^{\natural}_{\diamond} H}*id_{A^{\natural}_{\diamond} H})(a\o h)\\
 &&~~~~~\stackrel{}{=}(S_H(a_{1-1}\b\1(a_{2-1}\b\1(h_{1})))_{1}\rhd S_A(\a\2(a_{10})))(\b\1(S_H(a_{1-1}\b\1(a_{2-1}\\
 &&~~~~~~~~~\times \b\1(h_{1})))_{2})_{1}\rhd \a\2(a_{20}))\o \b\1(\b\1(S_H(a_{1-1}\b\1(a_{2-1}\b\1(h_{1})))_{2})_{2})h_{2}\\
 &&~\stackrel{(HA1)(HA2)}{=}(S_H(\b\1(a_{1-1}a_{2-1})\b\1(h_{1}))_{1}\rhd S_A(\a\2(a_{10})))(\b\1(S_H(\b\1(a_{1-1}a_{2-1})\\
 &&~~~~~~~~~\times \b\1(h_{1}))_{2})_{1}\rhd \a\2(a_{20}))\o \b\1(\b\1(S_H(\b\1(a_{1-1}a_{2-1})\b\1(h_{1}))_{2})_{2})h_{2}\\
 &&~~\stackrel{(HCMC1)}{=}(S_H(\b(a_{-1})\b\1(h_{1}))_{1}\rhd S_A(\a\2(a_{01})))(\b\1(S_H(\b(a_{-1})\\
 &&~~~~~~~~~\times \b\1(h_{1}))_{2})_{1}\rhd \a\2(a_{02}))\o \b\1(\b\1(S_H(\b(a_{-1})\b\1(h_{1}))_{2})_{2})h_{2}\\
 &&~\stackrel{(HC1)(HC2)}{=}(\b\1(S_H(\b(a_{-1})\b\1(h_{1}))_{11})\rhd S_A(\a\2(a_{01})))(\b\1(S_H(\b(a_{-1})\\
 &&~~~~~~~~~\times \b\1(h_{1}))_{12})\rhd \a\2(a_{02}))\o \b\1(S_H(\b(a_{-1})\b\1(h_{1}))_{2})h_{2}\\
 &&~\stackrel{(HC1)(HMA1)}{=}(\b(S_H(\b(a_{-1})\b\1(h_{1}))_{1})\rhd (S_A(\a\2(a_{01}))\a\2(a_{02}))\\
 &&~~~~~~~~~\o \b\1(S_H(\b(a_{-1})\b\1(h_{1}))_{2})h_{2}\\
 &&~~~\stackrel{(HA1)}{=}(\b(S_H(\b(a_{-1})\b\1(h_{1}))_{1})\rhd 1_A\v_A(a_{0}))\o \b\1(S_H(\b(a_{-1})\b\1(h_{1}))_{2})h_{2}\\
 &&~~\stackrel{(HCMC2)}{=}(\b(S_H(h_{1})_{1})\rhd 1_A\v_A(a))\o \b\1(S_H(h_{1})_{2})h_{2}\\
 &&~~~\stackrel{(HMA2)}{=}1_A\v_A(a)\o S_H(h_{1})h_{2}\\
 &&~~~~\stackrel{}{=}(1_{A}\o 1_{H})\v_A(a)\v_H(h)
 \end{eqnarray*}
 and
  \begin{eqnarray*}
 &&(id_{A^{\natural}_{\diamond} H}*S_{A^{\natural}_{\diamond} H})(a\o h)\\
 &&~~~~~\stackrel{}{=}a_1((a_{2-1}\b\1(h_{1}))_{1}\rhd \a\1(S_H(\a\1(a_{20})_{-1}\b\1(h_{2}))_{1}\rhd S_{A}(\a\2(\a\1(a_{20})_{0}))))\\
 &&~~~~~~~~~\o \b\1((a_{2-1}\b\1(h_{1}))_{2})\b\1(S_H(\a\1(a_{20})_{-1}\b\1(h_{2}))_{2})\\
 &&~\stackrel{(HM1)}{=}a_1((a_{2-1}\b\1(h_{1}))_{1}\rhd (\b\1(S_H(\a\1(a_{20})_{-1}\b\1(h_{2}))_{1})\rhd S_{A}(\a\3(\a\1(a_{20})_{0}))))\\
 &&~~~~~~~~~\o \b\1((a_{2-1}\b\1(h_{1}))_{2})\b\1(S_H(\a\1(a_{20})_{-1}\b\1(h_{2}))_{2})\\
 &&\stackrel{(HM2)(HA1)}{=}a_1(\b\1((a_{2-1}\b\1(h_{1}))_{1}S_H(\a\1(a_{20})_{-1}\b\1(h_{2}))_{1})\rhd S_{A}(\a\2(\a\1(a_{20})_{0})))\\
 &&~~~~~~~~~\o \b\1((a_{2-1}\b\1(h_{1}))_{2}S_H(\a\1(a_{20})_{-1}\b\1(h_{2}))_{2})\\
 &&~\stackrel{(HC1)}{=}a_1(\b\1((a_{2-1}\b\1(h_{1}))S_H(\a\1(a_{20})_{-1}\b\1(h_{2})))_{1}\rhd S_{A}(\a\2(\a\1(a_{20})_{0})))\\
 &&~~~~~~~~~\o \b\1((a_{2-1}\b\1(h_{1}))S_H(\a\1(a_{20})_{-1}\b\1(h_{2})))_{2}\\
 &&~\stackrel{(HCM1)}{=}a_1(\b\1((a_{2-1}\b\1(h_{1}))S_H(\b\1(a_{20-1})\b\1(h_{2})))_{1}\rhd S_{A}(\a\3(a_{200})))\\
 &&~~~~~~~~~\o \b\1((a_{2-1}\b\1(h_{1}))S_H(\b\1(a_{20-1})\b\1(h_{2})))_{2}\\
 &&~\stackrel{(HCM2)}{=}a_1(\b\1((\b\1(a_{2-11})\b\1(h_{1}))S_H(\b\1(a_{2-12})\b\1(h_{2})))_{1}\rhd S_{A}(\a\2(a_{20})))\\
 &&~~~~~~~~~\o \b\1((\b\1(a_{2-11})\b\1(h_{1}))S_H(\b\1(a_{2-12})\b\1(h_{2})))_{2}\\
 &&~\stackrel{(HC1)}{=}a_1((1_H\rhd S_{A}(\a\2(a_{20})))\v_H(a_{2-1})\o 1_H\v_H(h)\\
 &&~\stackrel{(HCM2)}{=}a_1(1_H\rhd S_{A}(\a\1(a_{2})))\o 1_H\v_H(h)\\
 &&~\stackrel{(HM2)}{=}a_1S_{A}(a_{2})\o 1_H\v_H(h)\\
 &&~~~~~\stackrel{}{=}(1_{A}\o 1_{H})\v_A(a)\v_H(h),
 \end{eqnarray*}
 while

  \begin{eqnarray*}
 &&S_{A^{\natural}_{\diamond} H}(\a(a)\o \b(h))\\
 &&~~~~\stackrel{}{=}(S_H(\a(a)_{-1}h)_1\rhd S_A(\a\2(\a(a)_0)))\o \b\1(S_H(\a(a)_{-1}h)_2)\\
 &&~\stackrel{(HCM1)}{=}(S_H(\b(a_{-1})h)_1\rhd S_A(\a\1(a_0)))\o \b\1(S_H(\b(a_{-1})h)_2)\\
 &&~~~~\stackrel{}{=}(\a\o \b)(S_{A^{\natural}_{\diamond} H}(a\o h)),
 \end{eqnarray*}
 finishing the proof.   \hfill $\square$
 \smallskip

  {\bf Corollary 3.7}(see \cite{MLY}) Let $(A,\a, S_A), (H,\b,S_H)$ be two Hom-Hopf algebras, and  $(A\natural H, \a\o \b)$ a smash product Hom-bialgebra.  Then $(A\natural H, \a\o \b, S_{A\natural H})$ is a Hom-Hopf algebra, where
  $$
  S_{A\natural H}(a\o h)=(S_H(h)_{1}\rhd \a\1(S_A(a)))\o \b\1(S_H(h)_{2}).
  $$

 {\bf Proof} Let the comodule action $\rho$ be trivial, i.e. $\rho(a)=1_H\o \a(a)$ in Theorem 3.6.         \hfill $\square$

 \smallskip

  {\bf Corollary 3.8}  Let $(C,\a,S_C), (H,\b,S_H)$ be two Hom-Hopf algebras, and $(C\diamond H, \a\o \b)$ a smash coproduct Hom-bialgebra.  Then $(C\diamond H, \a\o \b, S_{C\diamond H})$ is a Hom-Hopf algebra, where
  $$
  S_{C\diamond H}(c\o h)=S_C(\a\1(c_{(0)}))\o S_H(c_{(-1)}\b\1(h)).
  $$

 {\bf Proof} Let the module action $\rhd$ be trivial, i.e. $h\rhd c=\v_H(h)\a(c)$ in Theorem 3.6.         \hfill $\square$
 \smallskip

\section{Hom-Yetter-Drinfeld category}
\def\theequation{4. \arabic{equation}}
\setcounter{equation} {0} \hskip\parindent
 In this section, we give the definition of Hom-Yetter-Drinfeld module and also prove that the category $_H^H{\mathbb{YD}}$ of Hom-Yetter-Drinfeld modules is a pre-braided tensor category. Furthermore, we obtain that $(A^{\natural}_{\diamond} H,\a\o \b)$ is a Radford biproduct Hom-bialgebra if and only if  $(A,\a)$ is a Hom-bialgebra in the category $_H^H{\mathbb{YD}}$.

 \smallskip

 {\bf Definition 4.1} Let $(H, \b)$ be a Hom-bialgebra, $(M,\rhd_M, \a_M)$ a  left $(H,\b)$-module with action $\rhd_M: H\o M\lr M, h\o m\mapsto h\rhd_M m$ and $(M,\rho^M, \a_M)$ a left $(H,\b)$-comodule with coaction $\rho^M: M\lr H\o M, m\mapsto m_{-1}\o m_{0}$. Then we call $(M,\rhd_M, \rho^M,\a_M)$ a (left-left) Hom-Yetter-Drinfeld module over $(H,\b)$ if the following condition holds:
 $$
 (HYD)~~~h_1\b(m_{-1})\o (\b^3(h_2)\rhd_M m_0)=(\b^2(h_1)\rhd_M m)_{-1}h_2\o (\b^2(h_1)\rhd_M m)_{0},
 $$
 where $h\in H$ and $m\in M$.

 \smallskip

 {\bf Remarks}(1)The compatibility condition $(HYD)$ is different from the condition (2.1) in \cite[Definition 2.1]{MP1},  the condition (3.1) in \cite[Definition 3.1]{CZ} and  the condition (4.1) in \cite[Definition 4.1]{LS}.\newline
 \indent{\phantom{\bf Remarks}}(2) When $\b=id_H$,  the condition $(HYD)$ is exactly the condition $(YD)$.\newline
 \indent{\phantom{\bf Remarks}}(3) Let $(H, \b)$ be a Hom-bialgebra and $K$ a field. Then $(K, id_K)$ is a (left-left) Hom-Yetter-Drinfeld module over $(H,\b)$ with the module and comodule actions defined as follows: $H\o K\lr K, h\o k\mapsto \v(h)k$ and $K\lr H\o K, k\mapsto 1_H\o k$.\newline
 \indent{\phantom{\bf Remarks}}(4) When $(H,\b,S_H)$ is a Hom-Hopf algebra, then the condition $(HYD)$ is equivalent to
 $$
 (HYD)'~(\b^4(h)\rhd_M m)_{-1}\o (\b^4(h)\rhd_M m)_{0}=\b\2(h_{11}\b(m_{-1}))S_H(h_2)\o (\b^3(h_{12})\rhd_M m_0).
 $$
 \begin{eqnarray*}
 & {\bf Proof}&(\Longrightarrow)~~\b\2(h_{11}\b(m_{-1}))S(h_2)\o (\b^3(h_{12})\rhd m_0)\\
 &&~~~~~~~~~~~~\stackrel{(HYD)}{=}\b\2((\b^2(h_{11}\rhd m))_{-1}h_{12})S(h_{2})\o (\b^2(h_{11}\rhd m))_{0}\\
 &&~~~~~~\stackrel{(HA1)(HA2)}{=}\b\1((\b^2(h_{11}\rhd m))_{-1})(\b\2(h_{12})\b\1(S(h_{2})))\o (\b^2(h_{11}\rhd m))_{0}\\
 &&~~~~~~~~~\stackrel{(HC2)}{=}\b\1((\b^2(h_{1}\rhd m))_{-1})(\b\2(h_{21})\b\2(S(h_{22})))\o (\b^2(h_{1}\rhd m))_{0}\\
 &&~~~~~~~~~\stackrel{(HA1)}{=}\b\1((\b^2(h_{1}\rhd m))_{-1})(\b\2(h_{21}S(h_{22})))\o (\b^2(h_{1}\rhd m))_{0}\\
 &&~~~~~~\stackrel{(HA2)(HC2)}{=}(\b^4(h)\rhd m)_{-1}\o (\b^4(h)\rhd m)_{0}.\\
 &&(\Longleftarrow)~(\b^2(h_1)\rhd m)_{-1}h_2\o (\b^2(h_1)\rhd m)_{0}\\
 &&~~~~~~~~~\stackrel{(HYD)'}{=}(\b\2(\b\2(h_{1})_{11}\b(m_{-1}))S(\b\2(h_{1})_{2}))h_{2}\o (\b^3(\b\2(h_{1})_{12})\rhd m_{0})\\
 &&~~~~~~~~~\stackrel{(HC1)}{=}(\b\2(\b\2(h_{111})\b(m_{-1}))S(\b\2(h_{12})))h_{2}\o (\b(h_{112})\rhd m_{0})\\
 &&~~~~~\stackrel{(HC2)(HC1)}{=}(\b\2(\b\1(h_{11})\b(m_{-1}))S(\b\2(h_{21})))\b\1(h_{22})\o (\b^2(h_{12})\rhd m_{0})\\
 &&~~~~~\stackrel{(HA2)(HA1)}{=}(\b\1(\b\1(h_{11})\b(m_{-1}))(\b\2S(h_{21})h_{22})\o (\b^2(h_{12})\rhd m_{0})\\
 &&~~~~~~~~~~\stackrel{}{=}(\b\1(\b\1(h_{11})\b(m_{-1}))1_H\v_H(h_2)\o (\b^2(h_{12})\rhd m_{0})\\
 &&~\stackrel{(HC1)(HC2)(HA1)}{=}h_1\b(m_{-1})\o (\b^3(h_2)\rhd m_0).
 \end{eqnarray*}
 Here we use $\rhd$, $S$ instead of $\rhd_M$, $S_H$, respectively.             \hfill $\square$

  \smallskip

 {\bf Definition 4.2} Let $(H, \b)$ be a Hom-bialgebra. We denote by $_H^H{\mathbb{YD}}$ the category whose objects are Hom-Yetter-Drinfeld modules $(M,\rhd_M, \rho^M, \a_M)$ over $(H, \b)$; the morphisms in the category are morphisms of left $(H, \b)$-modules and left $(H, \b)$-comodules.

 \smallskip

 In the following, we give a solution of Hom-Yang-Baxter equation introduced and studied by Yau in \cite{Yau3,Yau7,Yau6}.

 {\bf Proposition 4.3} Let $(H, \b)$ be a Hom-bialgebra and $(M,\rhd_M, \rho^M, \a_M)$, $(N,\rhd_N, \rho^N, \a_N)$ $\in _H^H{\mathbb{YD}}$. Define the linear map
 $$
 \tau_{M,N}: M\o N\lr N\o M, ~~m\o n\mapsto \b^{3}(m_{-1})\rhd_N n \o m_{0},
 $$
 where $m\in M$ and $n\in N$. Then, we have $\tau_{M,N}\ci (\a_M\o \a_N)=(\a_N\o \a_M)\ci \tau_{M,N}$ and, if $(P,\rhd_P, \rho^P, \a_P) \in _H^H{\mathbb{YD}}$, the maps $\tau_{\underline{~~},~\underline{~~}}$ satisfy the Hom-Yang-Baxter equation:
 $$
 (\a_P\o \tau_{M,N})\ci (\tau_{M,P}\o \a_N)\ci (\a_M\o \tau_{N,P})=(\tau_{N,P}\o \a_M)\ci (\a_N\o \tau_{M,P})\ci (\tau_{M,N}\o \a_P).
 $$

 {\bf Proof} We only check that the second equality holds, and the first one is easy. For all $m\in M$, $n\in N$ and $p\in P$, we have
 \begin{eqnarray*}
 &&(\a_P\o \tau_{M,N})\ci (\tau_{M,P}\o \a_N)\ci (\a_M\o \tau_{N,P})(m\o n\o p)\\
 &&~~~~~\stackrel{}{=}(\b^{3}(\a_M(m)_{-1})\rhd_P (\b^{3}(n_{-1})\rhd_P p))\o \b^{3}(\a_M(m)_{0-1})\rhd_N \a_N(n_{0})\o \a_M(m)_{00}\\
 &&~~~\stackrel{(HM1)}{=}(\b^{4}(\a_M(m)_{-1})\rhd_P (\b^{4}(n_{-1})\rhd_P \a_P(p)))\o \b^{3}(\a_M(m)_{0-1})\rhd_N \a_N(n_{0})\\
 &&~~~~~~~~~~~~~\o \a_M(m)_{00}\\
 &&~~~\stackrel{(HCM1)}{=}(\b^{5}(m_{-1})\rhd_P (\b^{4}(n_{-1})\rhd_P \a_P(p)))\o \b^{4}(m_{0-1})\rhd_N \a_N(n_{0})\o \a_M(m_{00})\\
 &&~~~\stackrel{(HCM2)}{=}(\b^{4}(m_{-11})\rhd_P (\b^{4}(n_{-1})\rhd_P \a_P(p)))\o \b^{4}(m_{-12})\rhd_N \a_N(n_{0})\o \a_M^2(m_{0})\\
 &&~~~\stackrel{(HM2)}{=}((\b^{3}(m_{-11})\b^{4}(n_{-1}))\rhd_P \a_P^2(p))\o \b^{4}(m_{-12})\rhd_N \a_N(n_{0})\o \a_M^2(m_{0})\\
 &&~~~\stackrel{(HCM1)}{=}((\b^{3}(m_{-11}\a_N(n)_{-1}))\rhd_P \a_P^2(p))\o \b^{4}(m_{-12})\rhd_N \a_N(n)_{0}\o \a_M^2(m_{0})\\
 &&~~~\stackrel{(HA1)}{=}(\b^{2}(\b(m_{-11})\b(\a_N(n)_{-1})))\rhd_P \a_P^2(p))\o \b^{3}(\b(m_{-12}))\rhd_N \a_N(n)_{0}\o \a_M^2(m_{0})\\
 &&~~~\stackrel{(HC1)}{=}(\b^{2}(\b(m_{-1})_{1}\b(\a_N(n)_{-1})))\rhd_P \a_P^2(p))\o \b^{3}(\b(m_{-1})_{2})\rhd_N \a_N(n)_{0}\o \a_M^2(m_{0})\\
 &&~~~\stackrel{(HYD)}{=}(\b^{2}((\b^2(\b(m_{-1})_{1})\rhd_N \a_N(n))_{-1}\b(m_{-1})_{2})\rhd_P \a_P^2(p))\\
 &&~~~~~~~~~~~~~\o \b^2(\b(m_{-1})_{1})\rhd_N \a_N(n))_{0}\o \a_M^2(m_{0})\\
 &&~~~\stackrel{(HA1)(HC1)}{=}((\b^{2}((\b^3(m_{-11})\rhd_N \a_N(n))_{-1})\b^3(m_{-12}))\rhd_P \a_P^2(p))\\
 &&~~~~~~~~~~~~~\o (\b^3(m_{-11})\rhd_N \a_N(n))_{0}\o \a_M^2(m_{0})\\
 &&~~~\stackrel{(HCM2)}{=}((\b^{2}((\b^4(m_{-1})\rhd_N \a_N(n))_{-1})\b^3(m_{0-1}))\rhd_P \a_P^2(p))\\
 &&~~~~~~~~~~~~~\o (\b^4(m_{-1})\rhd_N \a_N(n))_{0}\o \a_M(m_{00})\\
 &&~~~\stackrel{(HM2)}{=}(\b^{3}((\b^4(m_{-1})\rhd_N \a_N(n))_{-1})\rhd_P (\b^3(m_{0-1})\rhd_P \a_P(p)))\\
 &&~~~~~~~~~~~~~\o (\b^4(m_{-1})\rhd_N \a_N(n))_{0}\o \a_M(m_{00})\\
 &&~~~\stackrel{(HM1)}{=}(\b^{3}(\a_N(\b^3(m_{-1})\rhd_N n)_{-1})\rhd_P (\b^3(m_{0-1})\rhd_P \a_P(p)))\\
 &&~~~~~~~~~~~~~\o \a_N(\b^3(m_{-1})\rhd_N n)_{0}\o \a_M(m_{00})\\
 &&~~~\stackrel{}{=}(\tau_{N,P}\o \a_M)\ci (\a_N\o \tau_{M,P})\ci (\tau_{M,N}\o \a_P)(m\o n\o p).
 \end{eqnarray*}   \hfill $\square$

 {\bf Lemma 4.4} Let $(H, \b)$ be a Hom-bialgebra and $(M,\rhd_M, \rho^M, \a_M)$, $(N,\rhd_N, \rho^N, \a_N)$ $\in _H^H{\mathbb{YD}}$. Define the linear maps
 $$
 \rhd_{M\o N}: H\o M\o N\lr M\o N, h\o m\o n\mapsto (h_1\rhd_M m)\o (h_2\rhd_N n),
 $$
 and
 $$
 \rho^{M\o N}: M\o N\lr H\o M\o N, m\o n\mapsto \b^{-2}(m_{-1}n_{-1})\o m_{0}\o n_{0},
 $$
 where $h\in H$, $m\in M$ and $n\in N$. Then $(M\o N, \rhd_{M\o N}, \rho^{M\o N}, \a_M\o \a_N)$ is a Hom-Yetter-Drinfeld module.

 {\bf Proof} It is easy to check that $(M\o N, \rhd_{M\o N}, \a_M\o \a_N)$ is an $(H, \b)$-Hom-module and $(M\o N, \rho^{M\o N}, \a_M\o \a_N)$ is an $(H, \b)$-Hom-comodule. While for $h\in H$, $m\in M$ and $n\in N$, we have
 \begin{eqnarray*}
 &&(\b^2(h_{1})\rhd_{M\o N} (m\o n))_{-1}h_{2}\o (\b^2(h_{1})\rhd_{M\o N} (m\o n))_{0}\\
 &&~~~~~\stackrel{}{=}((\b^2(h_{1})_{1}\rhd_{M} m)\o (\b^2(h_{1})_{2}\rhd_{N} n))_{-1}h_{2}\\
 &&~~~~~~~~~~~~~\o ((\b^2(h_{1})_{1}\rhd_{M} m)\o (\b^2(h_{1})_{2}\rhd_{N} n))_{0}\\
 &&~~~~~\stackrel{}{=}\b^{-2}(((\b^2(h_{1})_{1}\rhd_{M} m)_{-1}(\b^2(h_{1})_{2}\rhd_{N} n)_{-1})\b^2(h_{2}))\o (\b^2(h_{1})_{1}\rhd_{M} m)_{0}\\
 &&~~~~~~~~~~~~~\o (\b^2(h_{1})_{2}\rhd_{N} n)_{0}\\
 &&~~\stackrel{(HA1)(HA2)}{=}\b^{-2}(\b((\b^2(h_{11})\rhd_{M} m)_{-1})((\b^2(h_{12})\rhd_{N} n)_{-1}\b(h_{2})))\o (\b^2(h_{11})\rhd_{M} m)_{0}\\
 &&~~~~~~~~~~~~~\o (\b^2(h_{12})\rhd_{N} n)_{0}\\
 &&~~~~~\stackrel{(HC2)}{=}\b^{-2}(\b((\b^3(h_{1})\rhd_{M} m)_{-1})((\b^2(h_{21})\rhd_{N} n)_{-1}h_{22}))\o (\b^3(h_{1})\rhd_{M} m)_{0}\\
 &&~~~~~~~~~~~~~\o (\b^2(h_{21})\rhd_{N} n)_{0}\\
 &&~~~~~\stackrel{(HYD)}{=}\b^{-2}(\b((\b^3(h_{1})\rhd_{M} m)_{-1})(h_{21}\b(n_{-1})))\o (\b^3(h_{1})\rhd_{M} m)_{0}\\
 &&~~~~~~~~~~~~~\o (\b^3(h_{22})\rhd_{N} n_{0})\\
 &&~~~~~\stackrel{(HA2)}{=}\b^{-2}(((\b^3(h_{1})\rhd_{M} m)_{-1}h_{21})\b^2(n_{-1}))\o (\b^3(h_{1})\rhd_{M} m)_{0}\o (\b^3(h_{22})\rhd_{N} n_{0})\\
 &&~~~~~\stackrel{(HC2)}{=}\b^{-2}(((\b^2(h_{11})\rhd_{M} m)_{-1}h_{12})\b^2(n_{-1}))\o (\b^2(h_{11})\rhd_{M} m)_{0}\\
 &&~~~~~~~~~~~~~\o (\b^4(h_{2})\rhd_{N} n_{0})\\
 &&~~~~~\stackrel{(HYD)}{=}\b^{-2}((h_{11}\b(m_{-1}))\b^2(n_{-1}))\o (\b^3(h_{12})\rhd_{M} m_{0})\o (\b^4(h_{2})\rhd_{N} n_{0})\\
 &&~~~~~\stackrel{(HA1)}{=}(\b^{-2}(h_{11})\b\1(m_{-1}))n_{-1}\o (\b^3(h_{12})\rhd_{M} m_{0})\o (\b^4(h_{2})\rhd_{N} n_{0})\\
 &&~~~~~\stackrel{(HA2)}{=}\b^{-1}(h_{11})(\b\1(m_{-1})\b\1(n_{-1}))\o (\b^3(h_{12})\rhd_{M} m_{0})\o (\b^4(h_{2})\rhd_{N} n_{0})\\
 &&~~~~~\stackrel{(HC2)}{=}h_{1}(\b\1(m_{-1})\b\1(n_{-1}))\o (\b^3(h_{21})\rhd_{M} m_{0})\o (\b^3(h_{22})\rhd_{N} n_{0})\\
 &&~~\stackrel{(HC1)(HA1)}{=}h_{1}\b(\b^{-2}(m_{-1}n_{-1}))\o (\b^3(h_{2})_1\rhd_{M} m_{0})\o (\b^3(h_{2})_2\rhd_{N} n_{0})\\
 &&~~~~~\stackrel{}{=}h_{1}\b((m\o n)_{-1})\o (\b^3(h_{2})\rhd_{M\o N} (m\o n)_{0}),
 \end{eqnarray*}
 thus, the condition $(HYD)$ holds.  Therefore $(M\o N, \rhd_{M\o N}, \rho^{M\o N}, \a_M\o \a_N)$ is a Hom-Yetter-Drinfeld module.
   \hfill $\square$

 \smallskip

 {\bf Lemma 4.5} Let $(H, \b)$ be a Hom-bialgebra and $(M,\rhd_M, \rho^M, \a_M)$, $(N,\rhd_N, \rho^N, \a_N)$, $(P,\rhd_P, \rho^P, \a_P)$ $\in _H^H{\mathbb{YD}}$. With notation as
 above, define the linear map
 $$
 a_{M,N,P}: (M\o N)\o P\lr M\o (N\o P),~~(m\o n)\o p\mapsto \a_M^{-1}(m)\o (n \o \a_P(p)),
 $$
 where $m\in M$, $n\in N$ and $p\in P$. Then $a_{M,N,P}$ is an isomorphism of left $(H, \b)$-Hom-modules and left $(H, \b)$-Hom-comodules.

 {\bf Proof} Same to the proof of \cite[Proposition 3.2]{MP1}.
   \hfill $\square$

 \smallskip

 {\bf Lemma 4.6} Let $(H, \b)$ be a Hom-bialgebra and $(M,\rhd_M, \rho^M, \a_M)$, $(N,\rhd_N, \rho^N, \a_N)$ $\in _H^H{\mathbb{YD}}$. Define the linear map
 $$
 c_{M,N}: M\o N\lr N\o M,~~m\o n\mapsto (\b^2(m_{-1})\rhd_N \a_{N}^{-1}(n))\o \a_M^{-1}(m_{0}) ,
 $$
 where $m\in M$ and $n\in N$. Then $c_{M,N}$ is a morphism of left $(H, \b)$-Hom-modules and left $(H, \b)$-Hom-comodules.

 {\bf Proof} For all $h\in H$, $m\in M$ and $n\in N$, firstly,
 \begin{eqnarray*}
 &&(\a_N\o \a_M)\ci c_{M, N}(m\o n)\\
 &&~~~~~\stackrel{}{=}\a_N(\b^2(m_{-1})\rhd_N \a_N^{-1}(n))\o m_{0}\\
 &&~~~~~\stackrel{(HM1)}{=}(\b^3(m_{-1})\rhd_N  n)\o m_{0}\\
 &&~~~~~\stackrel{(HCM1)}{=}(\b^2(\a_M(m)_{-1})\rhd_N  \a_N^{-1}(\a_N(n))\o \a_M^{-1}(\a_M(m)_{0})\\
 &&~~~~~\stackrel{}{=}c_{M, N}\ci (\a_M\o \a_N)(m\o n);
 \end{eqnarray*}

 secondly,
 \begin{eqnarray*}
 &&c_{M, N}(h\rhd_{M\o N} (m\o n))=c_{M,N}((h_1\rhd_M m)\o (h_2\rhd_N n))\\
 &&~~~~~\stackrel{}{=}(\b^2((h_1\rhd_M m)_{-1})\rhd_N \a_N\1(h_2\rhd_N n))\o \a_M\1((h_1\rhd_M m)_{0})\\
 &&~~~~~\stackrel{(HM1)}{=}(\b^2((h_1\rhd_M m)_{-1})\rhd_N (\b\1(h_2)\rhd_N \a_N\1(n)))\o \a_M\1((h_1\rhd_M m)_{0})\\
 &&~~~~~\stackrel{(HM2)}{=}((\b((h_1\rhd_M m)_{-1})\b\1(h_2))\rhd_N n)\o \a_M\1((h_1\rhd_M m)_{0})\\
 &&~~~~~\stackrel{(HA1)}{=}(\b((h_1\rhd_M m)_{-1}\b^{-2}(h_2))\rhd_N n)\o \a_M\1((h_1\rhd_M m)_{0})\\
 &&~~~~~\stackrel{(HYD)}{=}(\b(\b^{-2}(h)_1\b(m_{-1}))\rhd_N n)\o \a_M\1(\b^3(\b^{-2}(h)_2)\rhd_M m_{0})\\
 &&~~~~~\stackrel{(HC1)}{=}(\b(\b^{-2}(h_1)\b(m_{-1}))\rhd_N n)\o \a_M\1(\b^3(\b^{-2}(h_2))\rhd_M m_{0})\\
 &&~~~~~\stackrel{}{=}((\b\1(h_1)\b^2(m_{-1}))\rhd_N n)\o \a_M\1(\b(h_2)\rhd_M m_{0})\\
 &&~~~~~\stackrel{(HM1)}{=}((\b\1(h_1)\b^2(m_{-1}))\rhd_N n)\o (h_2\rhd_M \a_M\1(m_{0}))\\
 &&~~~~~\stackrel{(HM2)}{=}(h_1\rhd_N (\b^2(m_{-1})\rhd_N \a_N\1(n)))\o (h_2\rhd_M \a_M\1(m_{0}))\\
 &&~~~~~\stackrel{}{=}h\rhd_{N\o M} ((\b^2(m_{-1})\rhd_N \a_N\1(n))\o \a_M\1(m_{0}))\\
 &&~~~~~\stackrel{}{=}h\rhd_{N\o M} c_{M,N}(m\o n);
 \end{eqnarray*}

 finally,
 \begin{eqnarray*}
 &&(\rho^{N\o M}\ci c_{M, N})(m\o n))\\
 &&~~~~~\stackrel{}{=}\b\2((\b^2(m_{-1})\rhd_N \a_N\1(n))_{-1}\a_M\1(m_{0})_{-1})\o (\b^2(m_{-1})\rhd_N \a_N\1(n))_{0}\\
 &&~~~~~~~~~~~~~\o \a_M\1(m_{0})_{0}\\
 &&~~~~~\stackrel{(HCM1)}{=}\b\2((\b^2(m_{-1})\rhd_N \a_N\1(n))_{-1}\b\1(m_{0-1}))\o (\b^2(m_{-1})\rhd_N \a_N\1(n))_{0}\\
 &&~~~~~~~~~~~~~\o \a_M\1(m_{00})\\
 &&~~~~~\stackrel{(HCM2)}{=}\b\2((\b(m_{-11})\rhd_N \a_N\1(n))_{-1}\b\1(m_{-12}))\o (\b(m_{-11})\rhd_N \a_N\1(n))_{0}\\
 &&~~~~~~~~~~~~~\o m_{0}\\
 &&~~~~~\stackrel{(HC1)}{=}\b\2((\b^2(\b\1(m_{-1})_{1})\rhd_N \a_N\1(n))_{-1}\b\1(m_{-1})_{2})\\
 &&~~~~~~~~~~~~~\o (\b^2(\b\1(m_{-1})_{1})\rhd_N \a_N\1(n))_{0}\o m_{0}\\
 &&~~~~~\stackrel{(HYD)}{=}\b\2(\b\1(m_{-1})_{1}\b(\a_N\1(n)_{-1}))\o (\b^3(\b\1(m_{-1})_{2})\rhd_N \a_N\1(n)_{0})\o m_{0}\\
 &&~~~~~\stackrel{(HC1)(HA1)}{=}\b\3(m_{-11})\b\1(\a_N\1(n)_{-1})\o (\b\2(m_{-12})\rhd_N \a_N\1(n)_{0})\o m_{0}\\
 &&~~~~~\stackrel{(HCM1)}{=}\b\3(m_{-11})\b\2(n_{-1})\o (\b\2(m_{-12})\rhd_N \a_N\1(n_{0}))\o m_{0}\\
 &&~~~~~\stackrel{(HCM2)}{=}\b\2(m_{-1})\b\2(n_{-1})\o (\b\2(m_{0-1})\rhd_N \a_N\1(n_{0}))\o \a_M\1(m_{00})\\
 &&~~~~~\stackrel{(HA1)}{=}\b\2(m_{-1}n_{-1})\o (\b\2(m_{0-1})\rhd_N \a_N\1(n_{0}))\o \a_M\1(m_{00})\\
 &&~~~~~\stackrel{}{=}(id\o c_{M,N})(\b\2(m_{-1}n_{-1})\o m_0\o n_0)\\
 &&~~~~~\stackrel{}{=}(id\o c_{M,N})\ci \rho^{M\o N}(m\o n).
 \end{eqnarray*}
Thus $c_{M,N}$ is a morphism of left $(H, \b)$-Hom-modules and left $(H, \b)$-Hom-comodules.                       \hfill $\square$

\smallskip

 {\bf Remarks} (1) The pre-braiding $(c_{M,N})$ differs from the one in \cite[Proposition 3.3]{MP1}.\newline
 \indent{\phantom{\bf Examples}} (2) If $(H, \b)$ is a Hom-Hopf algebra with bijective antipode $S$, then the pre-braiding $(c_{M,N})$ is invertible with 
 $$
 c_{M,N}^{-1}: N\o M\rightarrow M\o N, n\o m\mapsto \a_M^{-1}(m_{(0)})\o S^{-1}(\b^{2}(m_{(-1)}))\triangleright \a_N^{-1}(n).
 $$

 \smallskip

 {\bf Theorem 4.7} Let $(H, \b)$ is a Hom-Hopf algebra with bijective antipode $S$. Then the Hom-Yetter-Drinfeld category $_H^H{\mathbb{YD}}$ is a braided tensor category, with tensor product, associativity constraints, and braiding defined in Lemmas 4.4, 4.5 and 4.6, respectively, and the unit $I=(K, id_K)$.

 {\bf Proof} The proof of the pentagon axiom for $a_{M,N,P}$ is same to the proof of \cite[Theorem 3.4]{MP1}. Next we prove that the hexagonal relation for $c_{M,N}$.   Let $(M,\rhd_M, \rho^M, \a_M)$, $(N,\rhd_N, \rho^N, \a_N)$, $(P,\rhd_P, \rho^P, \a_P)$ $\in _H^H{\mathbb{YD}}$. Then for all $m\in M$, $n\in N$ and $p\in P$, we have
 \begin{eqnarray*}
 &&((id_N\o c_{M, P})\ci (a_{N,M,P})\ci (c_{M,N}\o id_P))((m\o n)\o p)\\
 &&~~~~\stackrel{}{=}\a_N\1(\b^2(m_{-1})\rhd_N \a_N\1(n))\o ((\b^2(\a_M\1(m_{0})_{-1})\rhd_P p)\o \a_M\1(\a_M\1(m_{0})_{0}))\\
 &&~\stackrel{(HCM1)}{=}\a_N\1(\b^2(m_{-1})\rhd_N \a_N\1(n))\o ((\b(m_{0-1})\rhd_P p)\o \a_M\2(m_{00}))\\
 &&~\stackrel{(HCM2)}{=}\a_N\1(\b(m_{-11})\rhd_N \a_N\1(n))\o ((\b(m_{-12})\rhd_P p)\o \a_M\1(m_{0}))\\
 &&~\stackrel{(HC1)}{=}\a_N\1(\b(m_{-1})_{1}\rhd_N \a_N\1(n))\o ((\b(m_{-1})_{2}\rhd_P p)\o \a_M\1(m_{0}))\\
 &&~\stackrel{(HCM1)}{=}\a_N\1(\b^2(\a_M\1(m)_{-1})_{1}\rhd_N \a_N\1(n))\o ((\b^2(\a_M\1(m)_{-1})_{2}\rhd_P p)\o \a_M\1(m)_{0})\\
 &&~~~~\stackrel{}{=}(a_{N,P,M}\ci c_{M,N\o P}\ci a_{M,N,P})((m\o n)\o p),
 \end{eqnarray*}
 and
 \begin{eqnarray*}
 &&((c_{M, P}\o id_N)\ci (a_{N,M,P}\1)\ci (id_M\o c_{N,P}))(m\o (n\o p))\\
 &&~~~~\stackrel{}{=}((\b^2(\a_M(m)_{-1})\rhd_P \a_P\1(\b^2(n_{-1})\rhd_P \a_P\1(p)))\o \a_M\1(\a_M(m)_{0}))\o \a_N\2(n_0)\\
 &&~~\stackrel{(HM1)}{=}((\b^2(\a_M(m)_{-1})\rhd_P (\b(n_{-1})\rhd_P \a_P\2(p)))\o \a_M\1(\a_M(m)_{0}))\o \a_N\2(n_0)\\
 &&~~\stackrel{(HM2)}{=}(((\b(\a_M(m)_{-1})\b(n_{-1}))\rhd_P \a_P\1(p))\o \a_M\1(\a_M(m)_{0}))\o \a_N\2(n_0)\\
 &&~\stackrel{(HM1)(HA1)}{=}(\a_P((\a_M(m)_{-1}n_{-1}))\rhd_P \a_P\2(p))\o \a_M\1(\a_M(m)_{0}))\o \a_N\2(n_0)\\
 &&~~~~\stackrel{}{=}(a_{P,M,N}\1\ci c_{M\o N,P}\ci a_{M,N,P}\1)(m\o (n\o p)),
 \end{eqnarray*}
 finishing the proof.                       \hfill $\square$

 \smallskip

 By Theorem 3.3, 3.6 and 4.7, we can get the main result in this paper.

 {\bf Theorem 4.8}  Let $(H, \b)$ is a Hom-Hopf algebra with bijective antipode $S$, $(A, \a)$ a left $(H, \b)$-module Hom-algebra and a left $(H, \b)$-comodule Hom-coalgebra satisfying $\b^2=id_H$. Then $(A^{\natural}_{\diamond} H, \mu_{A\natural H}, 1_A\o 1_H, \D_{A\diamond H}, \v_A\o \v_H, \a\o \b)$ is a Radford biproduct Hom-bialgebra if and only if  $(A,\a)$ is a bialgebra in the Hom-Yetter-Drinfeld category $_H^H{\mathbb{YD}}$.

 {\bf Proof}  It is obvious if we compare the conditions $(R4),(R5)$ in Theorem 3.3 and the condition $(HYD)$ in Definition 4.1, the definition of pre-braiding $c_{M,N}$ in Lemma 4.6, respectively.   \hfill $\square$

  {\bf Remarks} (1) If $\a=id_A$ and $\b=id_H$ in Theorem 4.8, then we can get the Majid's conclusion about the usual Radford biproduct and Yetter-Drinfeld category. \newline
 \indent{\phantom{\bf Remarks}} (2) $(A^{\natural}_{\diamond} H, \mu_{A\natural H}, 1_A\o 1_H, \D_{A\diamond H}, \v_A\o \v_H, \a\o \b, S_{A^{\natural}_{\diamond} H})$ is a Radford biproduct Hom-Hopf algebra if and only if  $(A,\a,S_A)$ is a Hopf algebra in the Hom-Yetter-Drinfeld category $_H^H{\mathbb{YD}}$.

 \smallskip

\section{Applications}
\def\theequation{5. \arabic{equation}}
\setcounter{equation} {0} \hskip\parindent
 In this section, we give some applications of the above sections.

 {\bf Example 5.1}  Let $K\mathbb{Z}_2=K\{1,a\}$ be Hopf group algebra (see \cite{Sw}). Then $(K\mathbb{Z}_2, id_{K\mathbb{Z}_2})$ is a Hom-Hopf algebra.

 Let $T_{2,-1}=K\{1, g, x, y|g^2=1, x^2=0,y=gx, gy=-gy=x\}$ be Taft's Hopf algebra (see \cite{MW1}), its coalgebra structure and antipode are given by
 $$
 \D(g)=g\o g,~\D(x) = x\o g+1\o x, ~\D(y)=y\o 1+g\o y;
 $$
 $$
 \v(g) = 1, \v(x) = 0, \v(y)=0;
 $$
 and
 $$
 S(g)=g,~S(x)=y,~S(y)=-x.
 $$

Define a linear map $\a$: $T_{2,-1}\lr T_{2,-1}$ by
 $$
 \a(1)=1,~\a(g)=g,~\a(x)=kx,~\a(y)=ky
 $$
 where $0\neq k\in K$. Then $\a$ is an automorphism of Hopf algebras.

 So we can get a Hom-Hopf algebra $H_{\a}=(T_{2,-1}, \a\ci \mu_{T_{2,-1}}, 1_{T_{2,-1}}, \D_{T_{2,-1}}\ci \a, \v_{T_{2,-1}}, \a)$ (see \cite{MS2}).

 \smallskip

 {\bf Lemma 5.1.1} With notations above.  Define module action $\rhd: K\mathbb{Z}_2\o H_{\a} \lr H_{\a}$ by
 $$
 1_{K\mathbb{Z}_2}\rhd 1_{H_{\a}}=1_{H_{\a}},~1_{K\mathbb{Z}_2}\rhd g=g,
 $$
  $$
 1_{K\mathbb{Z}_2}\rhd x=kx,~1_{K\mathbb{Z}_2}\rhd y=ky,
 $$
 $$
 a\rhd 1_{H_{\a}}=1_{H_{\a}},~a\rhd g=g,
 $$
 $$
 a\rhd x=kx,~a\rhd y=ky,
 $$
 Then by a routine computation we can get $(H_{\a}, \rhd, \a)$ is a $(K\mathbb{Z}_2, id_{K\mathbb{Z}_2})$-module Hom-algebra. Therefore,  $(H_{\a}\natural K\mathbb{Z}_2, \a\o id_{K\mathbb{Z}_2})$ is a smash product Hom-algebra.

 \smallskip

 {\bf Lemma 5.1.2} With notations above. Define comodule action $\rho: H_{\a} \lr K\mathbb{Z}_2\o H_{\a}$ by
 \begin{eqnarray*}
 &&\rho: H_{\a} \lr K\mathbb{Z}_2\o H_{\a}\\
 &&~~~~1_{H_{\a}}\mapsto 1_{K\mathbb{Z}_2}\o 1_{H_{\a}}\\
 &&~~~~~~g\mapsto 1_{K\mathbb{Z}_2}\o g\\
 &&~~~~~~x\mapsto ka\o x\\
 &&~~~~~~y\mapsto ka\o y.
 \end{eqnarray*}
 Then we can get $(H_{\a},\rho, \a)$ is a left $(K\mathbb{Z}_2, id_{K\mathbb{Z}_2})$-comodule Hom-coalgebra by a direct computation. Therefore, $(H_{\a}\natural K\mathbb{Z}_2, \a\o id_{K\mathbb{Z}_2})$ is a smash coproduct Hom-coalgebra.

 \smallskip

  By the above two lemmas and a direct computation, we have

  {\bf Theorem 5.1.3} With notations above. $({H_{\a}}^{\natural}_{\diamond} {K\mathbb{Z}_2}, \mu_{H_{\a}\natural {K\mathbb{Z}_2}},1_{H_{\a}}\o 1_{K\mathbb{Z}_2}, \D_{H_{\a}\diamond {K\mathbb{Z}_2}}, \v_{H_{\a}}\o \v_{K\mathbb{Z}_2}, \a\o id_{K\mathbb{Z}_2})$ is a Radford biproduct Hom-bialgebra. Furthermore, $({H_{\a}}^{\natural}_{\diamond} {K\mathbb{Z}_2}, \a\o id_{K\mathbb{Z}_2}, S_{{H_{\a}}^{\natural}_{\diamond} {K\mathbb{Z}_2}})$ is a Hom-Hopf algebra, where $S_{{H_{\a}}^{\natural}_{\diamond} {K\mathbb{Z}_2}}$ is defined by
 \begin{eqnarray*}
 &&S_{{H_{\a}}^{\natural}_{\diamond} {K\mathbb{Z}_2}}(1_{H_{\a}}\o 1_{K\mathbb{Z}_2})=1_{H_{\a}}\o 1_{K\mathbb{Z}_2};~~~S_{{H_{\a}}^{\natural}_{\diamond} {K\mathbb{Z}_2}}(1_{H_{\a}}\o a)=1_{H_{\a}}\o a\\
 &&~~~~~~~S_{{H_{\a}}^{\natural}_{\diamond} {K\mathbb{Z}_2}}(g\o 1_{K\mathbb{Z}_2})=g\o 1_{K\mathbb{Z}_2};~~~S_{{H_{\a}}^{\natural}_{\diamond} {K\mathbb{Z}_2}}(g\o a)=g\o a\\
 &&~~~~~~~~~~~~S_{{H_{\a}}^{\natural}_{\diamond} {K\mathbb{Z}_2}}(x\o 1_{K\mathbb{Z}_2})=y\o a;~~~S_{{H_{\a}}^{\natural}_{\diamond} {K\mathbb{Z}_2}}(x\o a)=y\o 1_{K\mathbb{Z}_2}\\
 &&~~~~~~~~~~S_{{H_{\a}}^{\natural}_{\diamond} {K\mathbb{Z}_2}}(y\o 1_{K\mathbb{Z}_2})=-x\o a;~~~S_{{H_{\a}}^{\natural}_{\diamond} {K\mathbb{Z}_2}}(y\o a)=-x\o 1_{K\mathbb{Z}_2}.
 \end{eqnarray*}

 \smallskip

 {\bf Example 5.2}  Let $K\mathbb{Z}_2=K\{1,a\}$ be Hopf group algebra (see \cite{Sw}). Then $(K\mathbb{Z}_2, id_{K\mathbb{Z}_2})$ is a Hom-Hopf algebra.

 Let $A=K\{1, z\}$ be a vector space. Define the multiplication $\mu_A$ by
 $$
 1z=z1=lz,~~z^2=0
 $$
 and the automorphism $\b: A\lr A$ by
 $$
 \b(1)=1,~~~\b(z)=lz
 $$
 where $0\neq l\in K$.  Then $(A,\b)$ is a Hom-algebra.

 Define the comultiplication $\D_A$ by
 $$
 \D_A(1)=1\o 1,~~\D_A(z)=lz\o 1+l1\o z,~~\hbox{and}~~\v_A(1)=1,~~\v_A(z)=0.
 $$
 Then $(A, \b)$ is a Hom-coalgebra.

 \smallskip

 {\bf Lemma 5.2.1} With notations above.  Define module action $\unrhd: K\mathbb{Z}_2\o A \lr A$ by
 $$
 1_{K\mathbb{Z}_2}\unrhd 1_A=1_{A},~1_{K\mathbb{Z}_2}\unrhd z=lz,
 $$
  $$
 a\unrhd 1_A=1_A,~a\unrhd z=-lz,
 $$
 Then by a routine computation we can get $(A, \unrhd, \b)$ is a $(K\mathbb{Z}_2, id_{K\mathbb{Z}_2})$-module Hom-algebra. Therefore,  $(A\natural K\mathbb{Z}_2, \b\o id_{K\mathbb{Z}_2})$ is a smash product Hom-algebra.

 \smallskip

 {\bf Lemma 5.2.2} With notations above. Define comodule action $\psi: A \lr K\mathbb{Z}_2\o A$ by
 \begin{eqnarray*}
 &&\psi: A \lr K\mathbb{Z}_2\o A\\
 &&~~~~1_{A}\mapsto 1_{K\mathbb{Z}_2}\o 1_{A}\\
 &&~~~~~~z\mapsto la\o z.
 \end{eqnarray*}
 Then we can get $(A,\psi, \b)$ is a left $(K\mathbb{Z}_2, id_{K\mathbb{Z}_2})$-comodule Hom-coalgebra by a direct computation. Therefore,
 $(A\natural K\mathbb{Z}_2, \b\o id_{K\mathbb{Z}_2})$ is a smash coproduct Hom-coalgebra.

 \smallskip

  By the above two lemmas and a direct computation, we have

  {\bf Theorem 5.2.3} With notations above. $({A}^{\natural}_{\diamond} {K\mathbb{Z}_2}, \mu_{A\natural {K\mathbb{Z}_2}},1_{A}\o 1_{K\mathbb{Z}_2}, \D_{A\diamond {K\mathbb{Z}_2}}, \v_{A}\o \v_{K\mathbb{Z}_2}, \b\o id_{K\mathbb{Z}_2})$ is a Radford biproduct Hom-bialgebra. Furthermore, $({A}^{\natural}_{\diamond} {K\mathbb{Z}_2}, \b\o id_{K\mathbb{Z}_2}, S_{{A}^{\natural}_{\diamond} {K\mathbb{Z}_2}})$ is a Hom-Hopf algebra, where $S_{{A}^{\natural}_{\diamond} {K\mathbb{Z}_2}}$ is defined by
 \begin{eqnarray*}
 &&S_{{A}^{\natural}_{\diamond} {K\mathbb{Z}_2}}(1_{A}\o 1_{K\mathbb{Z}_2})=1_{A}\o 1_{K\mathbb{Z}_2};~~~S_{{A}^{\natural}_{\diamond} {K\mathbb{Z}_2}}(1_{A}\o a)=1_{A}\o a\\
 &&~~~~~~~~S_{{A}^{\natural}_{\diamond} {K\mathbb{Z}_2}}(z\o 1_{K\mathbb{Z}_2})=z\o a;~~~S_{{A}^{\natural}_{\diamond} {K\mathbb{Z}_2}}(z\o a)=-z\o 1_{K\mathbb{Z}_2}.
 \end{eqnarray*}

 \smallskip

 {\bf Remark} If $\b=id_A$, i.e., $l=1$,  then Example 5.2 is same to the biproduct $B\star H$ (which is isomorphic to the Sweedler's Hopf algebra $T_{2,\omega}$) in \cite[Example 4.3]{MLZ}.

 \smallskip

 In the following, let us recall the definition of quasitriangular Hom-Hopf algebra in \cite{Yau3} or \cite{MLL}.

 A quasitriangular Hom-Hopf algebra is a octuple $(H,\mu,1_H, \D,$ $\v,S,\b, R)$ (abbr.$(H,\b,R)$) in which $(H,\mu,1_H,\D,\v,S,\b)$ is a Hom-Hopf algebra and  $R=R^1\o R^2 \in H\o H$, satisfying the following axioms (for all $h\in H$ and R=r):
 \begin{eqnarray*}
 &(QHA1)& \v(R^1)R^2=R^1\v(R^2)=1,\\
 &(QHA2)& {R^1}_{1}\o {R^1}_{2}\o \b(R^2)=\b(R^1)\o \b(r^1)\o R^2r^2,\\
 &(QHA3)& \b(R^1)\o {R^2}_{1}\o {R^2}_{2}=R^1r^1\o \b(r^2)\o \b(R^2),\\
 &(QHA4)& h_2R^1\o h_1R^2=R^1h_1\o R^2h_2,\\
 &(QHA5)& \b(R^1)\o \b(R^2)=R^1\o R^2.
 \end{eqnarray*}

 Let $(H,\b, S)$ be a Hom-Hopf algebra and $R=R^1\o R^2\in H\o H$. Define:
  $$
  \rho^H: H\lr H\o H~~~~~h\mapsto h_{-1}\o h_{0}=\b\3(R^2)\o R^1h.
  $$

  {\bf Proposition 5.3} Let $(H,\b, R)$ be a quasitriangular Hom-Hopf algebra. Then $(H, \b, \rho^H)$ is a left $(H,\b)$-comodule Hom-coalgebra and $(H, \mu_H,\rho^H,\b)$ is a Hom-Yetter-Drinfeld module.

  {\bf Proof} We compute as follows:
  \begin{eqnarray*}
  \b(h_{-1})\o \b(h_{0})
  &\stackrel{}{=}&\b(\b\3(R^2))\o \b(R^1h)\\
  &\stackrel{(HA1)}{=}&\b(\b\3(R^2))\o \b(R^1)\b(h)\\
  &\stackrel{(QHA5)}{=}&\b\3(R^2)\o R^1\b(h)=\b(h)_{-1}\o \b(h)_{0},
  \end{eqnarray*}
  so $(HCM1)$ holds.

  \begin{eqnarray*}
  h_{-11}\o h_{-12}\b(h_{0})
  &\stackrel{}{=}&\b\3(R^2)_1\o \b\3(R^2)_2\o \b(R^1h)\\
  &\stackrel{(HC1)(HA1)}{=}&\b\3({R^2}_1)\o \b\3({R^2}_2)\o \b(R^1)\b(h)\\
  &\stackrel{(QHA3)}{=}&\b\2({R^2})\o \b\2({r^2})\o (r^1R^1)\b(h)\\
  &\stackrel{(HA2)}{=}&\b\2({R^2})\o \b\2({r^2})\o \b(r^1)(R^1h)\\
  &\stackrel{(QHA5)}{=}&\b\2({R^2})\o \b\3({r^2})\o r^1(R^1h)\\
  &\stackrel{}{=}&\b(h_{-1})\o h_{0-1}\o h_{00},
  \end{eqnarray*}
  thus we get $(HCM2)$.

 \begin{eqnarray*}
  \b^2(h_{-1})\o h_{01}\o h_{02}
  &\stackrel{}{=}&\b\1(R^2)\o (R^1h)_1\o  (R^1h)_2\\
  &\stackrel{}{=}&\b\1(R^2)\o {R^1}_1h_1\o  {R^1}_2h_2\\
  &\stackrel{(QHA2)}{=}&\b\2(R^2r^2)\o \b({R^1})h_1\o  \b({r^1})h_2\\
  &\stackrel{(QHA5)(HA1)}{=}&\b\3(R^2)\b\3(r^2)\o {R^1}h_1\o  {r^1}h_2\\
  &\stackrel{}{=}&h_{1-1}h_{1-1}\o h_{10}\o h_{20},
  \end{eqnarray*}
  therefore we obtain $(HCMC1)$.

  $(HCMC2)$ can be checked by $(QHA1)$.

  Finally we verify that $(HYD)$ is satisfied.

 \begin{eqnarray*}
 (\b^2(h_{1})\rhd g)_{-1}h_{2}\o (\b^2(h_{1})\rhd g)_{0}
  &\stackrel{}{=}&\b\3(R^2)h_{2}\o R^1(\b^2(h_{1})g)\\
  &\stackrel{(HA2)}{=}&\b\3(R^2)h_{2}\o (\b\1(R^1)\b^2(h_{1}))\b(g)\\
  &\stackrel{(HA1)(HC1)}{=}&\b\3(R^2\b^3(h)_{2})\o \b\1(R^1\b^3(h)_{1})\b(g)\\
  &\stackrel{(QHA4)}{=}&\b\3(\b^3(h)_{1}R^2)\o \b\1(\b^3(h)_{2}R^1)\b(g)\\
  &\stackrel{(HA1)(HC1)}{=}&h_{1}\b\3(R^2)\o (\b^2(h_{2})\b\1(R^1))\b(g)\\
  &\stackrel{(HA2)}{=}&h_{1}\b\3(R^2)\o \b^3(h_{2})(\b\1(R^1)g)\\
  &\stackrel{(QHA5)}{=}&h_{1}\b\2(R^2)\o \b^3(h_{2})(R^1g)\\
  &\stackrel{}{=}&h_{1}\b(g_{-1})\o (\b^3(h_{2})\rhd g_{0}),
  \end{eqnarray*}
 finishing the proof.        \hfill $\square$

  \smallskip

 {\bf Proposition 5.4}  Let $(H,\b,S)$ be a Hom-Hopf algebra, with notations as above. If $(H, \b, \rho^H)$ is a left $(H,\b)$-comodule Hom-coalgebra and $(H, \mu_H,\rho^H,\b)$ is a Hom-Yetter-Drinfeld module, then $(H,\b, R)$ is a quasitriangular Hom-Hopf algebra.

 {\bf Proof} It is straightforward.     \hfill $\square$

  \smallskip

  By Proposition 5.3 and 5.4, we have:

  {\bf Theorem 5.5} With notations as above. $(H,\b, R)$ is a quasitriangular Hom-Hopf algebra if and only if $(H, \b, \rho^H)$ is a left $(H,\b)$-comodule Hom-coalgebra and $(H, \mu_H,\rho^H,\b)$ is a Hom-Yetter-Drinfeld module.

  \smallskip

  Dually, we have

  {\bf Theorem 5.6}  Let $(H,\b,S)$ be a Hom-Hopf algebra and $\s: H\o H\lr K$ a bilinear map. Define $\rhd_H: H\o H\lr H$ by
  $$
  h\o g\mapsto h\rhd_H g=\s(g_1, \b\3(h))g_2,
  $$
  where $h,g\in H$. Then $(H,\b, \s)$ is a cobraided Hom-Hopf algebra (see \cite{Yau4,MLY}) if and only if $(H, \b, \rhd_H)$ is a left $(H,\b)$-module Hom-algebra and $(H, \rhd_H,\D_H,\b)$ is a Hom-Yetter-Drinfeld module.

 {\bf Acknowledgments} This work was partially supported by the NNSF of China (No. 11101128), the Foundation of Henan Province (No. 14IRTSTHN023).

 \end{document}